\documentclass[12pt]{article}
\usepackage{amscd,amsfonts,amsmath,amssymb,amstext,amsthm}

\title{Hans Grauert: Mathematician Pur}
\author{Alan Huckleberry}
\date{\today}
\begin{document}
\maketitle
On April 25, 2008 on the occasion of the Gau\ss --Vorlesung in
Bonn Hans Grauert was presented the \emph{Ehrenmitgliedschaft}
of the Deutsche Mathematiker-Vereinigung (DMV). Given the opportuntiy
to describe even superficially the contributions of
this most distinguished colleague, I was pleased when asked 
to write this contribution for the Mitteilungen of the DMV. 
Only a few months after that article appeared\footnote{We wish
to thank the editors of the Mitteilungen for allowing us to
contribute the present article to the Notices. This is an
only slightly modified version of our article with the same title 
which appears in {\it Mitteilungen der Deutschen 
Mathematiker-Vereinigung, Band 16, Heft 2 (2008)}.}, 
on September 15, 2008 Grauert was awarded the prestigious 
Cantor-Medaille of the DMV at its annual meeting which took 
place this year in Erlangen. The inaugural Cantor-Medaille was
awarded in 1990 on the occasion of the 100th anniversary of
the founding of the DMV to Karl Stein who, like Grauert, devoted 
most of his scientific life to the subject of several complex 
variables.  The Medaille has been awarded on a regular basis 
since that time, being presented to J. Moser (1992), E. Heinz (1994), 
J. Tits (1996), V. Strassen (1999), Y. Manin (2002), F. Hirzebruch (2004) and
H. F\"ollmer (2006).  

\bigskip\noindent
It is \emph{presumptuous} for a mere mortal to even attempt to
write a Laudatio for Hans Grauert.  The easy part
is to construct a representative list of the various 
important stations of his life, and of his accomplishments
and honors, and to put them in a timeline. I have indeed
integrated a rough outline of these data in this article, but to me
this is just a small part of a story which I feel is very
important.  

\bigskip\noindent
Hans Grauert was born in Haren-Ems in 1930. At his retirement
festival in G\"ottingen he recalled how he struggled with
mathematics as a school boy until a teacher told him 
it was acceptable to think abstractly, he didn't 
necessarily need deal with numbers. 
No more than 15 years later he was introducing
spaces without points, just structure!

\bigskip\noindent
After beginning his studies in Sommersemester 1949 in Mainz, Grauert 
transferred to M\"unster, starting in Wintersemester 1949/50. 
There he was integrated into 
an exciting, energized mathematical atmosphere with friends
and teachers of all ages and experience. Among these was
Reinhold Remmert who would become his lifelong friend and main 
collaborator.  The mathematics guru was Karl Stein. Heinrich Behnke was well
connected to the outside mathematical world, in particular to
H. Hopf and H. Cartan, and had a very
good feeling for the important directions in complex analysis.

\bigskip\noindent
After a brief sojourn in Z\"urich, 
Grauert received his Dr. rer. nat. in M\"unster in 1954. Starting with his 
dissertation, Grauert contributed fundamental results which lie 
at the heart of a field of mathematics which
was in an infantile state when he started and was at a refined
and incredibly high level less than ten years later.
Let us now think back to the time when he began his studies!

\bigskip\noindent
There were indeed the deep, perhaps mysterious 
ideas of Oka on the table. Stein understood these in his own
way and was, for example, attempting to understand the role 
of topology in complex analysis, in particular for noncompact
spaces. Hirzebruch had received his doctorate in 1950 in M\"unster
and was on the path toward his fundamental book
\emph{Neue topologische Methoden in der algebraischen Geometrie}.
The power of the developing Cartan-Serre theory can not be underestimated.  
However, the foundations of what is now called
\emph{several complex variables} 
were simply not there!  

\bigskip\noindent
The works of Grauert and Remmert, together with the input
of Cartan and Serre (the postive
impact of the M\"unster-Paris connection is well-documented),
\emph{are} these foundations.  \emph{Komplexe R\"aume} and 
\emph{Bilder und Urbilder analytischer Garben} are two of numerous 
examples of their prolific joint work
which is basic for our subject. 

\bigskip\noindent
Perhaps also because they
have jointly written basic books on several complex variables,
\emph{Stein Theory} and \emph{Coherent Analytic Sheaves}, one 
might tend to overlook their different viewpoints. One can see 
that in this beginning phase Remmert is interested in analytic sets, 
their continuations, their properties with respect to holomorphic 
and meromorphic maps.

\bigskip\noindent
Grauert seems to be guided by problems involving complex
analytic objects on these sets. His \emph{Oka Principle},
which in terms which certainly hide the true depth of this work,
states that the category of holomorphic vector bundles
on a Stein space is the same as the category of 
topological bundles, is a perfect example.  The same
is true of his solution to the \emph{Levi Problem} where
he constructs holomorphic functions with given polar
data at the boundary of a domain by proving the
finite-dimensionality of a certain obstruction space
using a Fredholm theorem in a Fr\'echet context.  
His proof of the optimal version of \emph{Theorem B}
and his solution (with Docquier) of the Levi Problem 
for weakly pseudoconvex domains in Stein manifolds shows his 
deep understanding for approximation theorems of Runge
type.  

\bigskip\noindent
Shrinking coverings and understanding subtle properties
involving  restriction operators can be found at the top
of the list of Grauert's important methods.  The most complicated and perhaps
most famous of his results where such arguments appear
is his \emph{direct image theorem}.  Here one starts with
a proper holomorphic map $F:X\to Y$ between complex spaces where one
knows that the image $F(X)$ is an analytic subset of $Y$
(Remmert's Theorem).  Proving a theorem that is in a sense
in another universe, Grauert shows that direct images of
coherent sheaves on $X$ are coherent on $Y$. One
can not think of working in global complex geometry 
without the availability of this result! To obtain some feeling for
the order of magnitude of this work and for other interesting
information we recommend reading Remmert's article in
(\cite {R}).  

\bigskip\noindent
The last of the above mentioned works appeared in 1960, but
it is not at all clear to me when Grauert proved these 
theorems.  It seems that at a certain point he understood
\emph{everything} and it was just a matter of finding the
time and energy to write the papers. In any case he
choose the Oka Principle as his Habilitationsschrift and
around the time of completing his Hablitation, 
continuing the postwar tradition
which opened the world to numerous outstanding young German
mathematicians of that generation, Grauert left M\"unster
for the Institute of Advanced Study where he spent the
Wintersemester 1957/58. I know from direct discussions 
with others who were there at the time that the richness, 
depth and breadth of his ideas, which he 
presented both formally and informally, were nothing  
short of startling. 

\bigskip\noindent
In 1959 Grauert became Professor in G\"ottingen and
remained in this position until becoming Emeritus in 1996.
Once he told me he didn't want to be away from G\"ottingen 
for more than two weeks. But in fact he did travel widely. 
For example, most likely
due to the connection to Wilhelm Stoll, Grauert, Remmert
and Stein visited Notre Dame for extended periods. 
I know how important this was for that faculty and of 
course for me personally!

\bigskip\noindent
Grauert also invited distinguished foreign guests to G\"ottingen,
among them Aldo Andreotti. In the Winter Semester 1968/69 at 
Stanford I was introduced to Grauert's work in the lectures of Andreotti.  
Imagine being the only student in a course given by
the most wonderful of lecturers discusssing results
of his friend and coworker which are even in hindsight
some of the most beautiful in complex geometry. 
Their joint work was certainly one of the highlights: 
the Andreotti-Grauert theorems on finiteness and vanishing of cohomology 
for q-pseudoconvex manifolds, and their
jewel \emph{Algebraische K\"orper von automorphen Funktionen},
where they show how to use pseudoconcavity to prove the 
finite-dimensionality of spaces of automorphic forms.
However, what I remember most is Andreotti's explaining 
Grauert's elegant solution of the Levi problem and applications
to Kodaira-type vanishing and embedding theorems.

\bigskip\noindent
These last mentioned results are in a sense just snippets
of Grauert's remarkable paper \emph{\"Uber Modifikationen
und exzeptionelle analytische Mengen}. There, answers to fundamental
questions such as ``When can you blow down a variety?'' are
given. Concepts such as plurisubharmonicity, bundle curvature and signature of
intersection forms flow together.  A new, important criterion for projective
embeddings is proved.  After reading this work, I was sure
that this is the way mathematics should be!

\bigskip\noindent
On Grauert's research timeline we have now reached a point around 
1963. Of course the ideas kept coming!  There was a phase when
he was thinking about parameter spaces of complex analytic objects
(deformation theory). Here his two basic Inventiones papers
should be mentioned (\emph{\"Uber die Deformation isolierter
Singularit\"aten analytischer Mengen (1972)} and \emph{Der
Satz von Kuranishi f\"ur kompakte komplexe R\"aume (1974)}).
At the time when he was concentrating on vector bundles
(see for example his paper with M\"ulich, \emph{Vektorb\"undel
vom Rang 2 \"uber dem n-dimensionalen komplex-projektiven Raum (1975)}) , 
I remember a young mathematician asking him a general question
about what would nowadays be the most important direction 
of research in mathematics. Typifying how Grauert focused:
``Vector bundles on $\mathbb P_3$''! 

\bigskip\noindent
On the more analytic side there is the important work with  
his student Ramirez in the late 1960's and
then with Lieb on integral kernel representations.  The 
Grauert-Riemenschneider vanishing theorem (\emph{Verschwindungss\"atze f\"ur
analytische Kohomologiegruppen auf komplexen R\"aumen (1970)} 
can also be regarded as being at home in complex analysis.

\bigskip\noindent
More recently Grauert turned back to his old interests
in holomorphic and meromorphic equivalence relations. I remember 
he and Stein discussing these topics with great animation
just a few years before Stein's passing. His most
recent work in that direction in appeared in 1987. Finally,
one should not forget that \emph{hyperbolicity} has been in
been in the background for many years.  One sees this
in his work in 1965 with Reckziegel, his 1985 paper
with Ulrike Peternell, nee Grauert, and in his final paper
which was devoted to mathematics research 
\emph{Jetmetriken und hyperbolische Geometrie (1989)}.

\bigskip\noindent
At this point I could begin to be more precise about
the details of Grauert's work. However, I hope that the
above is sufficient for the interested bystander. For
those whose appetites might have been whetted, Grauert should
have the word: Please take a look at his collected works (\cite {G})
with its interesting annotations written by Yum-Tong Siu and 
Grauert himself.

\bigskip\noindent
As we all know, research is an extremely important part
of our academic lives, but there are other aspects which
must be emphasized. Here in Germany there is the classical notion of 
\emph{Akademischer Lehrer} which encompasses everything 
that a Professor should be.  Nowadays there seem to be
new interpretations being propagated.
\emph{Grants, research clusters, elite universities, etc} are
the buzzwords. However, we don't need new words to describe 
Grauert's contributions. Let me expand on this.

\bigskip\noindent
His work in administration of science must be commended, 
in particular his involvement with projects of the Deutsche
Forschungsgemeinschaft and his role on editorial boards, 
for example in bringing the Mathematische Annalen back to its historic 
high standards. Nevertheless, when I think of Grauert I think of him \emph{in}
the science not above it.  This includes his lectures, which 
may seem dry and minimal, sometimes even
formal, but you should listen very carefully. There are always 
deep ideas that should be followed! From the undergraduate
student in his Funktionentheorie course to the researcher
being advised in a private conversation, every listener 
should take every word seriously. The same is true of his 
vast written work. The reader must take the time to understand 
what is meant by every sentence!  This holds just as well for
his research monographs as it does for his textbooks. While
reading a Grauert-proof of Stokes' Theorem, you should keep
in mind that he has seriously thought about it!  These textbooks
range from the basic analysis sequence written with Fischer and Lieb
to the new version of Grauert-Fritzsche where even new
ideas in complex analysis are introduced.

\bigskip\noindent
Speaking of Grauert's minimality, I can't resist an anecdote.
Whenever he lectures he carries the \emph{Konzept} of the lecture with
him on a three by five card. He will most often start his lecture writing
\emph{Let $X$ be a complex space $\ldots$} on the
board and meticulously checking his Konzept to make sure
he got it right.  Given that he and Remmert originally
\emph{defined} the notion of a complex space, this
is a beautiful sight! Gossip has it that when giving a two 
semester course on several complex variables
he never changed the little piece of paper, but for the second
semester did in fact turn it over!!  

\bigskip\noindent
I have no idea how many students (Diplom, Staatsexam, Promotion)
have done their work with Grauert. In any case it is a large
number! We who are working in areas near their works see
the strong positive influence of the master teacher, and
I know that Grauert is proud of them all.  

\bigskip\noindent
A researcher of the highest quality, a teacher at all levels  
with relevant fundamental new ideas always in the background, 
an author with a style where every word has a meaning,
an important participant in and leader of academic societies, 
a cultured intellectual in the sense
of Humboldt, and a very kind gentleman, Grauert personifies the true 
notion of \emph{Akademischer Lehrer}.  


\begin {thebibliography} {XXX}
\bibitem [G] {G}
Grauert, H.: Selected papers I. and II., Springer-Verlag, 1994
%
\bibitem [R] {R}
Remmert, R.: Intelligencer, vol. 17, No2, 1995, 4-11
\end {thebibliography}
\end{document}